\def\eqdef{\buildrel \text def \over =}

\renewcommand{\phi}{\varphi}

\newcommand{\su}{\subseteq}

\newcommand{\lng}{\langle}
\newcommand{\rng}{\rangle}

\newcommand{\N}{\mathbb N}

\newcommand{\pr}{\operatorname{Pr}}

\documentclass[12pt,reqno]{amsart}
 \usepackage{amstext}
 \usepackage{amssymb}

\newtheorem{lemma}{Lemma}
\newtheorem{theorem}[lemma]{Theorem}

 \newtheorem{definition}[lemma]{Definition}
  \newtheorem{corollary}[lemma]{Corollary}
 
 \newtheorem{problem}[lemma]{Problem}
 \newtheorem{claim}[lemma]{Claim}

\begin{document}
\title[regressive Ramsey numbers]{Regressive Ramsey numbers are 
Ackermannian} 
\author{Menachem Kojman}
\address{Department of Mathematics and Computer Science\\
Ben Gurion University of the Negev\\
Beer Sheva, Israel}
\address{Department of Mathematical Sciences\\
Carnegie-Mellon University\\
Pittsburgh, PA}

\email{kojman@cs.bgu.ac.il}
\thanks{The first author was partially supported by  NSF grant 
No. DMS-9622579}

\author{Saharon Shelah}
\address{Institute of Mathematics\\
The Hebrew University of Jerusalem\\
Jerusalem 91904, Israel} 
\email{shelah@math.huji.ac.il} 
\thanks{The second author was partially supported by the Binational 
Science Foundation.  Number 649 in list of publications}

\begin{abstract} We give an elementary proof of the fact that 
regressive Ramsey numbers are Ackermannian. This fact  was first proved by 
Kanamori and McAloon  with mathematical logic techniques.
\end{abstract}

\maketitle

\begin{quote} 
\textit{\small
Nous vivons encore sous le r\`{e}gne de la logique, voil\`{a}, bien 
entendu, \`{a} quoi je voulais en venir. Mais les proc\'{e}d\'{e}s 
logiques, de nos jours, ne s'appliquent plus qu'\`{a} la r\'{e}solution 
de probl\`{e}mes d'int\'{e}r\^{e}t secondaire. \cite[1924, p. 
13]{manifesto}}
 \end{quote}

\section{Introduction}

\begin{definition}
\begin{enumerate}
\item	let $A$ be a set of natural numbers. A coloring $c:[A]^e\to \N$ 
	of unordered $e$-tuples from $A$ is {\em regressive} if $c(x)<\min x$ 
	for all $x\in [A]^e$.
	
\item A subset $B\su A$ is {\em min-homogeneous} for a coloring $c$ of 
$[A]^e$ if for all $x\in [A]^e$ the color $c(x)$ depends only on 
$\min x$.
\end{enumerate}
\end{definition}

\begin{theorem}[Kanamori and McAloon]

\begin{enumerate}
	\item For every $k$ and $e$ there exists $N$ such that for every 
	regressive  coloring of $e$-tuples from $\{1,2,\dots,N\}$ there exists a 
	min-homogeneous subset of size $k$.  
	
	\item The statement in (1) 
	cannot be proved from the axioms of Peano Arithmetic (although it 
	can be \emph{phrased} in the language of PA) 
	
	\item Let $\nu(k)$ be the 
	least $N$ which satisfies 1 for $e=2$.  The function $\nu$ 
	eventually dominates every primitive recursive function.
\end{enumerate}
\end{theorem}

Part (3) of Kanamori and McAloon's result \cite{KM} was proved with 
methods from mathematical logic.  We present below an elementary proof of 
3.

\section{The lower bound}
For every function $f:\N\to \N$  and  $n$, $f^{(n)}$ is
defined  by
$f^{(0)}(x)=x$ and $f^{(n+1)}(x)=f(f^{(n)}(x))$ for all $x\in \N$.

 Define a sequence of (strictly
increasing) integer  functions $f_i:\N\to \N$ for   $i\ge 1$ as follows:

	\begin{align}\label{defEq}
	f_1(n)&=n+1  \\
	f_{i+1}(n)&=f_i^{(\lfloor \sqrt{n} /2\rfloor)}(n)
	 \end{align}

Fix an integer $k>2$. 
Define a sequence of semi-metrics $\lng d_i:i\in \N\rng$ on 
$\{n:n \ge 4k^2\}$ by putting, for $m,n\ge4k^2$, 

\begin{gather}
d_i(m,n)=|\{l\in \N : m< f_i^{(l)}(4k^2)\le n\}|
\end{gather}

 For $n>m\ge4k^2$ let $I(m,n)$ be the greatest $i$ for which $d_i(m,n)$ 
is positive, and  $d(m,n)=d_{I(m,n)}(m,n)$.

\begin{claim} \label{sqrt}  For all $n\ge m\ge 4k^2$,  $d(m,n) \le 
\sqrt{m}/2$.
\end{claim}

\begin{proof}

  Let $i=I(m,n)$.  Since $d_{i+1}(m,n)=0$, there exists $t$ 
 and $l$ such that $t=f_{i+1}^{(l)}(4k^2)\le m \le 
 n<f_{i+1}^{(l+1)}(4k^2)=f_{i+1}(t)$.  But $f_{i+1}(t)=f_{i}^{(\lfloor 
 \sqrt{t}/2\rfloor)}(t)$ and therefore $\sqrt{t}/2\ge d_i(t,f_{i+1}(t))\ge 
 d(m,n)$.
\end{proof}

Let us fix the following (standard) pairing function $\pr$ on 
$\N^2$ 

\[\label{pairing} \pr(m,n)=\binom{m+n+1}{2}+n
\]

$\pr$ is a bijection between $[\N]^2$ and $\N$ and is monotone in each 
variable.  Observe that if $m,n\le l$ then $\pr(m,n)<4l^2$ for all 
$l>2$.

Define a pair coloring $c$ on $\{n:n\ge 4k^2\}$ as follows: 

\begin{gather}
\label{coloring}
c(\{m,n\})=\pr(I(m,n),d(m,n)) 
\end{gather}

\begin{claim} \label{noMinHom} 
For every $i\in\N$, every sequence $x_0<x_1<\dots<x_i$ that satisfies 
$d_i(x_0,x_i)=0$ is not min-homogeneous for $c$.
\end{claim}

\begin{proof}
The claim is proved by induction on $i$.  If $i=1$ then there are no 
$x_0<x_1$ with $d_1(x_0,x_1)=0$ at all.  Suppose 
to the contrary that $i>1$, that $x_0<x_1<\dots< x_i$ form a 
min-homogeneous sequence with respect to $c$ and that $d_i(x_0,x_i)=0$.  
Necessarily, $I(x_0,x_i)=j<i$. By min-homogeneity, 
$I(x_0,x_1)=j$ as well, and $d_j(x_0,x_i)=d_j(x_0,x_1)$.  Hence, 
$\{x_1,x_2,\dots x_i\}$ is min-homogeneous with $d_j(x_1,x_i)=0$ --- 
contrary to the induction hypothesis.
\end{proof}

\begin{claim}\label{regressive} The coloring 
$c$  is regressive on the interval $[4k^2,f_k(4k^2))$.
\end{claim}

\begin{proof} Clearly, $d_{k+1}(m,n)=0$ for $4k^2\le 
m<n<f_k(4k^2)$ and therefore $I(m,n)< k\le\sqrt{m}/2$.  From Claim \ref{sqrt} 
we know that $d(m,n)\le \sqrt{m}/2$.  Thus, 
$c(\{m,n\})
\le \pr(\lfloor\sqrt{m}/2\rfloor,\lfloor\sqrt{m}/2\rfloor)$, which is $<m$, since 
$\sqrt{m}>2$.
\end{proof}

We show that $f_k(4k^2)$ grows eventually faster than every primitive 
recursive function by comparing the functions $f_i$ with the usual 
approximations of Ackermann's function. It is well known that every 
primitive recursive function is dominated by some approximation of 
Ackermann's function (see, e.g. \cite{Calude}).

Let $A_i(n)$ be defined as follows: 

\begin{align}
A_1(n)&=n+1\\  
A_{i+1}(n)&=A_i^{(n)}(n)
\end{align}

 The $A_i$-s are the usual approximations to 
Ackermann's function, which is defined by $Ack(n)=A_n(n)$.

\begin{claim}
    \begin{enumerate}
       \item for all $n\ge 16$ and $i\ge 7$,
        \begin{enumerate} \item $16n^2\le f_i(n)$;
                          \item $f_i(16n^2)\le f_i^{(2)}(n)$.
        \end{enumerate}
       \item $f_i(n)\le f_{i+6}^{(2)}(n)$ for all $i\ge 1$ and $n\ge 16$. 
    \end{enumerate}
\end{claim}

\begin{proof}
Inequality (a) is verified directly.

Inequality (b)  follows from (a) by substituting $f_i(n)$ for
$16n^2$ in $f_i(16n^2)$, since $f_i$ is increasing.

We prove 2  by induction on $i$. For $i=1$ it holds that
$n+1<f_7^{(2)}(n)$ for all $n\ge 16$ by (a). 

Suppose the inequality holds for $i$ and all $n\ge 16$, and let $n\ge 16$ be given.
Since
$A_i(n)\le f_{i+6}^{(2)}(n)$ for all $n\ge 16$, it follows by monotonicity of
$A_i$ that
$A_i^{(n)}(n)\le f_{i+6}^{(2n)}(n)$. The latter term is smaller than
$f_{i+6}^{(2n)}(16n^2)$ by monotonicity, which equals $f_{i+7}(16n^2)$ by  (2).
Inequality (b)  implies that $A_i^{(n)}(n)\le f_{i+7}^{(2)}(n)$. Finally,
$A_i^{(n)}(n)=A_{i+1}(n)$ by (6).
\end{proof}

 \begin{claim}
 For all $i\ge 7$ and $n\ge 16$ it holds that $A_i(n)\le f_{i+7}(n)$.
 \end{claim}
 
 \begin{proof}
 By 2 in the previous claim, $A_i(n)\le f_{i+6}^{(2)}(n)$ for $n\ge 16$. If $n\ge
16$, then $\sqrt{n}/2\ge 2$ and hence, by (2),
$f_{i+7}(n)\eqdef f_{i+6}^{\lfloor(\sqrt{n}/2\rfloor)}\ge f_{i+6}^{(2)}(n)$.
 \end{proof}

\begin{corollary}
The function $\nu(k)$ eventually dominates every primitive recursive function.
\end{corollary}

\section{Discussion}

\subsection{Other Ramsey numbers}

Paris and Harrington \cite{PH} published in 1976 the first finite 
Ramsey-type statement that was shown to be independent over Peano 
Arithmetic.  Soon after the discovery of the Paris-Harrington result, 
Erd\H os and Mills studied the Ramsey-Paris-Harrington numbers in 
\cite{EM}.  Denoting by $R^e_c(k)$ the Ramsey-Paris-Harrington number 
for exponent $e$ and $c$ many colors, Erd\H os and Mills showed that 
$R^2_2(k)$ is double exponential in $k$ and that $R^2_c(k)$ is 
Ackermannian as a function of $k$ and $c$.  In the same paper, several 
small Ramsey-Paris-Harrington numbers were computed.  Later Mills tightened 
the double exponential upper bound for $R^2_2(k)$ in \cite{Mills}.

Canonical Ramsey numbers for pair colorings were treated in \cite{LR} 
 and were also found to be double exponential.

The second author showed that van der Waerden numbers are primitive 
recursive, refuting the conjecture that they were Ackermannian, in 
\cite{Sh:329} (see also \cite{Noga}).  

We remark that an upper bound for regressive Ramsey numbers for pairs 
is $R^3_2(k)$ --- the Ramsey-Paris-Harrington number for 
\emph{triples}.  Let $N$ be large enough and suppose that $c$ is 
regressive on $\{1,2,\dots,N-1\}$.  Color a triple $x<y<z$ red if 
$c(x,y)=c(x,z)$ and blue otherwise.  Find a homogeneous set $A$ of 
size at least $k$ and so that $|A|>\min A + 1$.  The homogeneous color 
on $A$ cannot be blue for $k>5$, and therefore $A$ is min-homogeneous 
for $c$.

\subsection{Problems}
   The following two problems about regressive Ramsey numbers remain 
   open:
 
 \begin{problem}
 \begin{enumerate}
 \item
  Find a concrete upper bound for regressive Ramsey 
 numbers.
 \item Compute  small regressive Ramsey numbers
 \end{enumerate}
 \end{problem}
 

\end{document}